\documentclass[english,11pt]{article}
\usepackage{amsfonts}
\usepackage[latin1]{inputenc}
\usepackage{amsmath,amsthm,amssymb}
\usepackage{graphicx}
\usepackage{color}

\def\neweq#1{\begin{equation}\label{#1}}
\def\endeq{\end{equation}}

\newtheorem{theorem}{Theorem}[section]
\newtheorem{lemma}{Lemma}[section]
\newtheorem{definition}{Definition}[section]

\newtheorem{remark}{Remark}[section]
\input{tcilatex}
\begin{document}

\title{\textbf{EXISTENCE AND ASYMPTOTIC BEHAVIOR OF SOLUTION TO A SINGULAR
ELLIPTIC PROBLEM} to appear in Surveys in Mathematics and its Applications }
\author{Dragos-Patru Covei\break \\
{\small \ Constantin Brancusi University of \ Tg-Jiu and West University of
Timisoara, Romania}\\
{\small E-mail: c\texttt{ovdra@yahoo.com}}}
\date{}
\maketitle

\begin{abstract}
In this paper we obtain existence results for the positive solution of a
singular elliptic boundary value problem. To prove the main results we use
comparison arguments and the method of sub-super solutions combined with a
procedure which truncates the singularity.
\end{abstract}

\baselineskip16pt \renewcommand{\theequation}{\arabic{section}.%
\arabic{equation}} \catcode`@=11 \@addtoreset{equation}{section} \catcode%
`@=12



\section{ Introduction}

This paper contains contribution of a technical nature to the study of
positive solutions of the equations 
\begin{equation}
-\Delta u+c(x)u^{-1}\left\vert \nabla u\right\vert ^{2}=a(x)\text{ for }x\in 
\mathbb{R}^{N}\text{, }u>0\text{ in }\mathbb{R}^{N}\text{, }u(x)\rightarrow 0%
\text{ as }\left\vert x\right\vert \rightarrow \infty  \label{A1}
\end{equation}%
where $N>2$, $a:\mathbb{R}^{N}\rightarrow \mathbb{R}$ is a function
satisfying the following conditions

AC1)\quad $a,c\in C_{loc}^{0,\alpha }(\mathbb{R}^{N})$ \ for some $\alpha
\in (0,1);$

AC2)\quad $a(x)>0,c(x)>0$ for all $x\in \mathbb{R}^{N};$

A3)\quad for $\varphi (r)=\max_{|x|=r}a(x)$ we have%
\begin{equation*}
\int_{0}^{\infty }r\varphi (r)dr<\infty .
\end{equation*}

Problems like (\ref{A1}) has been intensively studied. Our study is
motivated by the works of Shu \cite{WEN}, ~Arcoya, Carmona, Leonori,
Aparicio, Orsina and Petitta \cite{DAR}, Arcoya, Barile and Aparicio \cite%
{DAR2} where the existence, non-existence and uniqueness of solution for the
problem like (\ref{A1}) are solved.

In this article we present a new argument in the study of the problem (\ref%
{A1}) more simple that used in \cite{DAR}, \cite{DAR2}, \cite{WEN} and where
the problem is considered just in the case when $\Omega \subset \mathbb{R}%
^{N}$ is a bounded domain with smooth boundary.

The above equation contains different quantities, such as: singular
nonlinear term (like $u^{-1}$), convection nonlinearity (denoted by $%
\left\vert \nabla u\right\vert ^{2}$), as well as potentials ($c$ and $a$).
The principal difficulty in the treatment of (\ref{A1}) is due to the
singular character of the equation combined with the nonlinear gradient term.

The importance of the problem\textit{\ }(\ref{A1}) is given considering \
the well know problem%
\begin{equation}
\begin{array}{ll}
\Delta u=a(x)h(u)\text{, }u>0\text{ \textit{in} \ }\Omega ,\text{ } & \text{%
\ }u(x)=\infty \text{ \ \textit{as} \ }x\rightarrow \partial \Omega ,%
\end{array}
\label{bib}
\end{equation}%
\ because we can easily deduce the following two remarks:

\begin{remark}
\label{11}When $h(u)=e^{u}$, by a transformation of the form $w=e^{-u}$ the
problem (\ref{bib}) becomes 
\begin{equation}
-\Delta w+\frac{\left\vert \nabla w\right\vert ^{2}}{w}=a(x)\text{, }w>0%
\text{ }in\text{ }\Omega \text{, }w\left( x\right) \rightarrow 0\text{ }as%
\text{ \ }x\rightarrow \partial \Omega \text{,}  \label{1d}
\end{equation}%
but this is the problem (\ref{A1}) when $c(x)=1$.
\end{remark}

\begin{remark}
\label{12}For $h(u)=u^{\delta }$ $(\delta >1)$ and $w=C[u]^{-C^{-1}},(C:=1/(%
\delta -1))$ in (\ref{bib}) we have%
\begin{equation}
-\Delta w+\delta C\frac{\left\vert \nabla w\right\vert ^{2}}{w}=a(x)\text{, }%
w>0\text{, }in\text{ }\Omega \text{, }w\rightarrow 0\text{ }as\text{ \ }%
x\rightarrow \partial \Omega \text{,}  \label{2d}
\end{equation}%
which is the problem (\ref{A1}) when $c(x)=\delta C$.
\end{remark}

This finish the motivation of our work.

The main results of the article are:

\begin{theorem}
\label{idea1} If $\Omega \subset \mathbb{R}^{N}$ is a bounded domain with
boundary $\partial \Omega $ of class $C^{2,\alpha }$ for some $\alpha \in
(0,1)$ and $a,$ $c\in C^{0,\alpha }(\overline{\Omega }),$ $a(x)>0,$ $c(x)>0$
for any $x\in \overline{\Omega },$ then the problem%
\begin{equation}
\text{ }-\Delta u+c(x)u^{-1}\left\vert \nabla u\right\vert ^{2}=a(x)\text{
in\ \ }\Omega \text{, }u_{\left\vert \partial \Omega \right. }=0,
\label{gr1}
\end{equation}%
has at least a positive solution $u\in C(\overline{\Omega })\cap C^{2,\alpha
}(\Omega )$.
\end{theorem}

In the next result we establish sufficient condition for the existence of
solution to the problem (\ref{A1}) in the case when $\Omega =\mathbb{R}^{N}$.

\begin{theorem}
\label{1.1.1}We suppose that hypotheses AC1), AC2), A3) are satisfied. Then,
the problem (\ref{A1}) has a $C_{loc}^{2,\alpha }(\mathbb{R}^{N})$ positive
solution vanishing at infinity. If, in addition,%
\begin{equation}
\lim_{\left\vert x\right\vert \rightarrow \infty }\left\vert x\right\vert
^{\mu }\varphi (\left\vert x\right\vert )<\infty ,  \label{as11}
\end{equation}%
for some $\mu \in (2,N)$, then%
\begin{equation}
u(x)=O(\left\vert x\right\vert ^{2-\mu })\text{ as }\left\vert x\right\vert
\rightarrow \infty \text{.}  \label{as22}
\end{equation}
\end{theorem}

To prove the existence of such a solution to (\ref{A1}) we establish some
preliminary results.

\section{Preliminary results}

Since we apply sub and super solution method due to Amann \cite{AMANNC3}, we
recall the following definition of sub and super solution which are our main
tools in the proof of the solvability of problem (\ref{A1}).

For $f_{1}(x,\eta ,\xi ):\overline{\Omega }\times \mathbb{R}\times \mathbb{R}%
^{N}\rightarrow \mathbb{R}$ and $g_{1}:\partial \Omega \rightarrow \mathbb{R}
$, Amann introduce the following definitions:

\begin{definition}
A function $\underline{u}\in C^{2,\alpha }(\overline{\Omega })$ is called a
sub solution for the problem 
\begin{equation}
-\Delta u=f_{1}(x,u,\nabla u)\text{ in }\Omega ,\text{ }u=g\text{ on }%
\partial \Omega \text{,}  \label{amann}
\end{equation}%
if 
\begin{equation*}
-\Delta \underline{u}\leq f_{1}(x,\underline{u},\nabla \underline{u})\text{
in }\Omega ,\text{ }\underline{u}=g\text{ on }\partial \Omega .
\end{equation*}
\end{definition}

\begin{definition}
A function $\overline{u}\in C^{2,\alpha }(\overline{\Omega })$ is called a
super solution of the problem (\ref{amann}) if 
\begin{equation*}
-\Delta \overline{u}\geq f_{1}(x,\overline{u},\nabla \overline{u})\text{ in }%
\Omega ,\text{ }\overline{u}=g\text{ on }\partial \Omega .
\end{equation*}
\end{definition}

One of the important results from \cite{AMANNC3} is:

\begin{lemma}
\label{aman}Let $\Omega $ be a bounded domain from $\,\mathbb{R}^{N}$, with
boundary $\partial \Omega $ of class $C^{2,\alpha }$ for some $\alpha \in
(0,1)$, $g\in C^{2,\alpha }(\partial \Omega )$ and $f_{1}$ be a continuous
function with the property that $\partial f_{1}/\partial \eta $, $\partial
f_{1}/\partial \xi ^{i}$, $i=\overline{1,N}$ \ exists and are continuous on $%
\overline{\Omega }\times \mathbb{R}^{N+1}$ and such that

AM1) $f_{1}(\cdot ,\eta ,\xi )\in C^{\alpha }(\overline{\Omega })$,
uniformly for $(\eta ,\xi )$ in bounded subsets of $\mathbb{R\times R}^{N}$;

AM2)there exists a function $f_{2}:$\ $\mathbb{R}_{+}\rightarrow \mathbb{R}%
_{+}:=[0,\infty )$ such that%
\begin{equation}
\left\vert f_{1}(x,\eta ,\xi )\right\vert \leq f_{2}(\rho )(1+\left\vert \xi
\right\vert ^{2}),  \label{quadrically}
\end{equation}%
for every $\rho \geq 0$ and $(x,\eta ,\xi )\in \overline{\Omega }\times
\lbrack -\rho ,\rho ]\times \mathbb{R}^{N}$.

Under these assumption, if the problem (\ref{amann}) has a sub solution $%
\underline{u}$ and a super solution $\overline{u}$ such that $\underline{u}%
(x)\leq \overline{u}(x),$ $\forall x\in \overline{\Omega }$ then there
exists at least a function $u(x)\in C^{2+\alpha }(\overline{\Omega })$ which
satisfies $\underline{u}(x)\leq u(x)\leq \overline{u}(x)$ for all $x\in 
\overline{\Omega }$ and satisfying (\ref{amann}) pointwise.\bigskip\ More
precisely, there exist a minimal solution $\overset{\thicksim }{u}(x)\in
\lbrack \underline{u}(x),\overline{u}(x)]$ and a maximal solution $\overset{%
\approx }{u}(x)\in \lbrack \underline{u}(x),\overline{u}(x)]$, in the sense
that every solution $u(x)\in $ $[\underline{u}(x),\overline{u}(x)]$
satisfies $\overset{\thicksim }{u}(x)\leq u(x)\leq $ $\overset{\approx }{u}%
(x)$.
\end{lemma}

We will need the following variant of the maximum principle:

\begin{lemma}
\label{maximum}\bigskip\ Assume that $\Omega $ is a bounded open set in $%
\mathbb{R}^{N}$. If $u:\overline{%
\Omega
}\rightarrow \mathbb{R}$ is a smooth function such that%
\begin{equation*}
\left\{ 
\begin{array}{ll}
-\Delta u\geq 0 & in\text{ }\Omega , \\ 
u\geq 0 & on\text{ }\partial \Omega ,%
\end{array}%
\right.
\end{equation*}%
then $u\geq 0$ in $\Omega .$
\end{lemma}

This finishes the auxiliary results. Now we prove the announced Theorems.

\section{Proof of the Theorem \protect\ref{idea1}}

In the following will we use similarly argument that were used by Crandall,
Rabinowitz and Tartar \cite{CRTC3}, Noussair \cite{NS2} and the author \cite%
{CDC5}.

Let $\varepsilon \in (0,1)$. The existence will be established by solving
the approximate problems%
\begin{equation}
\left\{ 
\begin{array}{rr}
-\Delta u+c(x)u^{-1}\left\vert \nabla u\right\vert ^{2}=a(x), & in\text{\ \ }%
\Omega ,\text{ }u>\varepsilon \text{ }in\text{\ \ }\Omega , \\ 
u=\varepsilon , & \text{\textit{on} }\partial \Omega .%
\end{array}%
\right.  \label{grgr}
\end{equation}

For this, let $\varphi _{1}$ be the first positive eigenfunction
corresponding to the first eigenvalue $\lambda _{1}$ of the problem 
\begin{equation}
-\Delta u(x)=\lambda u(x),\text{ \textit{in} \ }\Omega ,\text{ }%
u_{\left\vert \partial \Omega \right. }(x)=0.  \label{eigenvalue1}
\end{equation}%
It is well known \ that $\varphi _{1}\in C^{2+\alpha }(\overline{\Omega })$.
We note by $m_{2}:=\min_{x\in \overline{\Omega }}a(x)$ and $%
M_{1}:=\max_{x\in \overline{\Omega }}c(x)$ to prove that the function $%
\underline{u}(x)=\sigma _{1}\varphi _{1}^{2}+\varepsilon ,$ where 
\begin{equation}
0<\sigma _{1}\leq \min \left\{ \frac{m_{2}}{2\lambda _{1}\max_{x\in 
\overline{\Omega }}\varphi _{1}^{2}+4M_{1}\max_{x\in \overline{\Omega }%
}\left\vert \nabla \varphi _{1}\right\vert ^{2}},1\right\}  \label{small1}
\end{equation}%
is a sub solution of (\ref{grgr}) in the sense of Lemma \ref{aman}. Indeed,
by (\ref{small1}) we have%
\begin{equation*}
\begin{array}{l}
-\Delta \underline{u}+c(x)\underline{u}^{-1}\left\vert \nabla \underline{u}%
\right\vert ^{2}-a(x)\leq -\Delta \underline{u}+M_{1}\underline{u}%
^{-1}\left\vert \nabla \underline{u}\right\vert ^{2}-m_{2} \\ 
\leq -2\sigma _{1}\varphi _{1}\Delta \varphi _{1}-2\sigma _{1}\left\vert
\nabla \varphi _{1}\right\vert ^{2}+4M_{1}\sigma _{1}\left\vert \nabla
\varphi _{1}\right\vert ^{2}-m_{2} \\ 
=2\sigma _{1}\lambda _{1}\varphi _{1}^{2}-2\sigma _{1}\left\vert \nabla
\varphi _{1}\right\vert ^{2}+4M_{1}\sigma _{1}\left\vert \nabla \varphi
_{1}\right\vert ^{2}-m_{2} \\ 
\leq 2\sigma _{1}\lambda _{1}\varphi _{1}^{2}+4M_{1}\sigma _{1}\left\vert
\nabla \varphi _{1}\right\vert ^{2}-m_{2}\leq 0.%
\end{array}%
\end{equation*}%
In the next step we prove the existence of a super solution to the problem (%
\ref{grgr}). For this, let $v\in C^{2+\alpha }(\overline{\Omega })$ be the
unique solution of the problem 
\begin{equation}
-\Delta y=a(x)\text{ \textit{in} \ }\Omega ,\text{ }y(x)=0\text{ for }x\in
\partial \Omega .  \label{Step11}
\end{equation}%
We observe that, $\overline{u}=v+\varepsilon \in C^{2+\alpha }(\overline{%
\Omega })$, fulfils%
\begin{equation*}
-\Delta \overline{u}(x)+c(x)\overline{u}^{-1}(x)\left\vert \nabla \overline{u%
}(x)\right\vert ^{2}=a(x)+c(x)\overline{u}^{-1}(x)\left\vert \nabla 
\overline{u}(x)\right\vert ^{2}\geq a(x)\text{ for }x\in \Omega .
\end{equation*}%
Clearly, $\overline{u}$ is a super solution to (\ref{grgr}). Now, since%
\begin{equation}
\left\{ 
\begin{array}{llll}
-\Delta \lbrack \overline{u}-\underline{u}] & \geq & a(x)+c(x)\underline{u}%
^{-1}\left\vert \nabla \underline{u}\right\vert ^{2}-a(x)\geq 0, & 
\mbox{ in
}\text{ \ }\Omega , \\ 
\text{ \ \ \ \ \ }\overline{u}-\underline{u} & = & 0, & \mbox{ on }\text{ }%
\partial \Omega ,%
\end{array}%
\right.  \label{max}
\end{equation}%
follows from the maximum principle, \textit{Lemma \ref{maximum},} that $%
\underline{u}(x)\leq \overline{u}(x)$, $x\in \overline{\Omega }$.

We have obtained a sub solution $\underline{u}\in C^{2,\alpha }(\overline{%
\Omega })$ and a super solution $\overline{u}\in C^{2,\alpha }(\overline{%
\Omega })$ for the problem (\ref{grgr}) such that $\underline{u}\leq 
\overline{u}$ in $\overline{\Omega }$ with the property from \textit{Lemma %
\ref{aman}}. Then, there exists $u_{\varepsilon }\in C^{2,\alpha }(\overline{%
\Omega })$ such that 
\begin{equation}
\begin{array}{cccccc}
\underline{u}(x) & \leq & u_{\varepsilon }(x) & \leq & \overline{u}(x), & 
\text{ }x\in \overline{\Omega }.%
\end{array}
\label{ineq1}
\end{equation}%
and satisfying (pointwisely) the problem (\ref{grgr}).

The relation (\ref{ineq1}) shows that $u>0$ in $\Omega $. We remark that $%
\underline{u}=\sigma _{1}v^{2}+\varepsilon ,$ where $\sigma _{1}$ is a
positive constant such that 
\begin{equation}
0<\sigma _{1}\leq \min \left\{ \frac{m_{2}}{\max_{x\in \overline{\Omega }%
}[2v+4M_{1}\left\vert \nabla v\right\vert ^{2}]},1\right\} ,  \label{small2}
\end{equation}%
is again a sub solution of (\ref{grgr}) with the same property from \textit{%
Lemma \ref{aman}}.

In this time we have obtained a function $u_{\varepsilon }\in C^{2,\alpha }(%
\overline{\Omega })$ that satisfies pointwisely the equivalently form of (%
\ref{grgr}): 
\begin{equation}
\left\{ 
\begin{array}{ll}
-\Delta u+c(x)\left( u+\varepsilon \right) ^{-1}\left\vert \nabla
u\right\vert ^{2}=a(x), & \text{\textit{in}\ }\Omega , \\ 
u>0, & \text{\textit{in}\ }\Omega , \\ 
u=0, & \text{\textit{on} }\partial \Omega .%
\end{array}%
\right.  \label{unicl}
\end{equation}%
Moreover $u_{\varepsilon }\in C^{2,\alpha }(\overline{\Omega })$ is unique.
Indeed, assume that the problem (\ref{unicl}) has more that one solution and
let $v_{\varepsilon }$ the second solution. Let us show that $u_{\varepsilon
}\leq v_{\varepsilon }$ or, equivalently, $u_{\varepsilon }\left( x\right)
+\varepsilon \leq v_{\varepsilon }\left( x\right) +\varepsilon $ for any $%
x\in \overline{\Omega }$. Assume the contrary. Set%
\begin{equation*}
\alpha (x):=\frac{u_{\varepsilon }\left( x\right) +\varepsilon }{%
v_{\varepsilon }\left( x\right) +\varepsilon }-1.
\end{equation*}%
Since we have $\left. \left[ \alpha \left( x\right) \right] \right\vert
_{\partial \Omega }=0$ we deduce that $\max_{\overline{\Omega }}\alpha
\left( x\right) $, exists and is positive. At that point, say $x_{0}$, we
have $\nabla \alpha (x_{0})=0$ \textit{and} $\Delta \alpha (x_{0})\leq 0$,
which implies%
\begin{equation}
\displaystyle\Big(-\left( v_{\varepsilon }+\varepsilon \right) \Delta
u_{\varepsilon }+\left( u_{\varepsilon }+\varepsilon \right) \Delta
v_{\varepsilon }\Big)(x_{0})\geq 0,  \label{max1}
\end{equation}%
and 
\begin{equation}
\displaystyle\frac{\left\vert \nabla u_{\varepsilon }(x_{0})\right\vert ^{2}%
}{\left( u_{\varepsilon }(x_{0})+\varepsilon \right) ^{2}}=\frac{\left\vert
\nabla v_{\varepsilon }\right\vert ^{2}}{\left( v_{\varepsilon
}(x_{0})+\varepsilon \right) ^{2}}.  \label{max2}
\end{equation}%
By $\left( \ref{max1}\right) $ and $\left( \ref{max2}\right) $ we have 
\begin{equation}
\displaystyle\frac{a\left( x_{0}\right) }{u_{\varepsilon
}(x_{0})+\varepsilon }-\frac{a\left( x_{0}\right) }{v_{\varepsilon
}(x_{0})+\varepsilon }+c(x_{0})\left( \frac{\left( v_{\varepsilon
}+\varepsilon \right) ^{-1}\left\vert \nabla v_{\varepsilon }\right\vert ^{2}%
}{v_{\varepsilon }+\varepsilon }-\frac{\left( u_{\varepsilon }+\varepsilon
\right) ^{-1}\left\vert \nabla u\right\vert ^{2}}{u_{\varepsilon
}+\varepsilon }\right) (x_{0})\geq 0,  \label{cont1}
\end{equation}%
or, equivalently%
\begin{equation}
\displaystyle a\left( x_{0}\right) \frac{v_{\varepsilon
}(x_{0})-u_{\varepsilon }(x_{0})}{\left( u_{\varepsilon }(x_{0})+\varepsilon
\right) \left( v_{\varepsilon }(x_{0})+\varepsilon \right) }\geq 0.
\label{cont3}
\end{equation}%
which is a contradiction with $u_{\varepsilon }(x_{0})>v_{\varepsilon
}(x_{0})$. So $u_{\varepsilon }(x)\leq v_{\varepsilon }(x)$ in $\overline{%
\Omega }$. A similar argument can be made to produce $v_{\varepsilon
}(x)\leq u_{\varepsilon }(x)$ forcing $u_{\varepsilon }(x)=v_{\varepsilon
}(x)$.

We will show that, for any smooth bounded subdomain $\Omega ^{\prime }$ of $%
\mathbb{R}^{N}$ there exists a constant $C_{4}>0$ such that%
\begin{equation}
\left\Vert u_{\varepsilon }\right\Vert _{C^{2,\alpha }(\overline{\Omega }%
^{\prime })}\leq C_{4}.  \label{claim}
\end{equation}%
For any bounded $C^{2,\alpha }$-smooth domain $\Omega ^{\prime }\subset $\ $%
\mathbb{R}^{N}$, take $\ \Omega _{1}$, $\Omega _{2}$ and $\Omega _{3}$ with $%
C^{2,\alpha }$-smooth boundaries, \ such that $\Omega ^{\prime }\subset
\subset \Omega _{1}\subset \subset \Omega _{2}\subset \subset \Omega
_{3}\subset \subset \Omega $. Note that%
\begin{equation}
u_{\varepsilon }(x)\geq \underline{u}\left( x\right) >0\text{, }\forall x\in
\Omega _{i},\text{ }i=\overline{1,3}\text{.}  \label{r1}
\end{equation}%
Let $h_{\varepsilon }(x)=a(x)-c(x)\left( u_{\varepsilon }\left( x\right)
+\varepsilon \right) ^{-1}\left\vert \nabla u_{\varepsilon }\left( x\right)
\right\vert ^{2},x\in \overline{\Omega _{3}}$. Following, we use $C_{i=%
\overline{1,4}}$, to denote positive constants which are independent of $%
\varepsilon $.

Since $-\Delta u_{\varepsilon }(x)=h_{\varepsilon }(x)$, $x\in \overline{%
\Omega _{3}},$ we see by the interior gradient estimate theorem of
Ladyzenskaya and Ural'tseva \cite[Theorem 3.1, p. 266]{LUC3} that there
exists a positive constant $C_{1}$ independent of $\varepsilon $\ such that 
\begin{equation}
\max_{x\in \overline{\Omega }_{2}}\left\Vert \nabla u_{\varepsilon }\left(
x\right) \right\Vert \leq C_{1}\max_{x\in \overline{\Omega }%
_{3}}u_{\varepsilon }\left( x\right) .  \label{CRT}
\end{equation}%
Using (\ref{ineq1}) and (\ref{CRT}) we obtain that $\left\Vert \nabla
u_{\varepsilon }\right\Vert $ is uniformly bounded on $\overline{\Omega }%
_{2} $. This final result, the property of $a$ and $c$ shows that $%
\left\vert h_{\varepsilon }\right\vert $ is uniformly bounded on $\overline{%
\Omega }_{2} $ and so $h_{\varepsilon }\in L^{p}(\Omega _{2})$ for any $p>1$.

Since $-\Delta u_{\varepsilon }(x)=h_{\varepsilon }(x)$ for $x\in \Omega
_{2} $, we see from\textit{\ \cite{CDC5},} that there exists a positive
constant $C_{2}$\ independent of $\varepsilon $ such that 
\begin{equation*}
\left\Vert u_{\varepsilon }\right\Vert _{W^{2,p}(\Omega _{1})}\leq
C_{2}(\left\Vert h_{\varepsilon }(x)\right\Vert _{L^{p}(\Omega
_{2})}+\left\Vert u_{\varepsilon }\right\Vert _{L^{p}(\Omega _{2})}),
\end{equation*}%
i.e. $\left\Vert u_{\varepsilon }\right\Vert _{W^{2,p}(\Omega _{1})}$ is
uniformly bounded.

Choose $p$ such that $p>N$ and $p>N\left( 1-\alpha \right) ^{-1}$. Then by
Sobolev's imbedding theorem, it follows that $\left\Vert u_{\varepsilon
}\right\Vert _{C^{1,\alpha }\left( \overline{\Omega }_{1}\right) }$ is
uniformly bounded by a constant independent of $\varepsilon $.

Moreover, this say that $h_{\varepsilon }\in C^{0,\alpha }(\overline{\Omega }%
_{1})$ and $\left\Vert h_{\varepsilon }\right\Vert _{C^{0,\alpha }(\overline{%
\Omega }_{1})},$ is uniformly bounded. Using this and the interior Schauder
estimates (see \textit{\cite{CDC5,GTC3})}, for solutions of elliptic
equations (\ref{dd101}) we have that there exists a positive constant $C_{3}$
independent of $\varepsilon $ with the property%
\begin{equation}
\left\Vert u_{\varepsilon }\right\Vert _{C^{2,\alpha }(\overline{\Omega }%
^{\prime })}\leq C_{3}\left( \left\Vert h_{\varepsilon }\right\Vert
_{C^{0,\alpha }(\overline{\Omega }_{1})}+\sup_{\overline{\Omega }%
_{1}}u_{\varepsilon }\right) .  \label{bound1}
\end{equation}%
Because $\left\Vert h_{\varepsilon }\right\Vert _{C^{0,\alpha }(\overline{%
\Omega }_{1})}$ is uniformly bounded, we see from (\ref{bound1}) that%
\begin{equation}
\left\Vert u_{\varepsilon }\right\Vert _{C^{2,\alpha }\left( \overline{%
\Omega }^{\prime }\right) }\leq C_{4}.  \label{2a}
\end{equation}%
Thus (\ref{claim}) is proved.

Set $\varepsilon :=1/n$ and $u_{\varepsilon }:=u^{n}$. Since the sequence $%
u^{n}$ is bounded in $C^{2,\alpha }\left( \overline{\Omega }^{\prime
}\right) $ for any bounded domain $\Omega ^{\prime }\subset \subset $\ $%
\Omega $ by (\ref{2a}), using the Ascoli-Arzela theorem and the standard
diagonal process, we can find a subsequence of $u^{n}$, denote again by $%
u^{n}$ and a function $u\in C^{2}\left( \overline{\Omega }^{\prime }\right) $
\ such that $\left\Vert u^{n}-u\right\Vert _{C^{2}\left( \overline{\Omega }%
^{\prime }\right) }\rightarrow 0$ for $n\rightarrow \infty $. In particular 
\begin{equation*}
\Delta u^{n}\ \text{\textit{respectively} }a(x)-c(x)(u^{n}(x)+1/n)^{-1}\left%
\vert \nabla u^{n}(x)\right\vert ^{2}\ 
\end{equation*}%
\ converge for $n\rightarrow \infty $ in $\overline{\Omega }^{\prime }$ to 
\begin{equation*}
\Delta u\text{ \textit{respectively} }a(x)-c(x)u(x)^{-1}\left\vert \nabla
u(x)\right\vert ^{2}.
\end{equation*}%
It follows that $u$ is a solution of%
\begin{equation}
-\Delta u=a(x)-c(x)u^{-1}(x)\left\vert \nabla u(x)\right\vert ^{2},\text{ 
\textit{in} }\overline{\Omega }^{\prime },  \label{probprob}
\end{equation}%
of class $C^{2}(\overline{\Omega }^{\prime }),$ and hence of class $%
C^{2,\alpha }(\overline{\Omega }^{\prime })$ by a standard regularity
arguments based on Schauder estimates.

Since $\Omega ^{\prime }$ is arbitrary, we also see that $u\in
C_{{}}^{2,\alpha }(\Omega )$. We have obtained $u^{n}\overset{n\rightarrow
\infty }{\rightarrow }u$ (pointwisely) in $C_{{}}^{2,\alpha }(\Omega )$.

For $\varepsilon :=1/n\overset{n\rightarrow \infty }{\rightarrow }0$ in (\ref%
{ineq1}) we have 
\begin{equation}
\begin{array}{cccccc}
\underline{u}_{2}(x):=\sigma _{1}\varphi _{1}^{2} & \leq & u(x) & \leq & 
\overline{u}^{2}(x):=v(x), & \text{ }x\in \overline{\Omega }.%
\end{array}
\label{finfin}
\end{equation}

Moreover, by (\ref{probprob}) and (\ref{finfin}), we obtain%
\begin{equation*}
-\Delta u=a(x)-c(x)u^{-1}\left\vert \nabla u\right\vert ^{2}\text{ \textit{%
a.e. in} }\Omega \text{, }u>0\text{ \textit{in} }\Omega \text{, }%
u_{\left\vert \partial \Omega \right. }=0.
\end{equation*}%
Thus $u\in C(\overline{\Omega })\cap C_{{}}^{2,\alpha }(\Omega )$ is the
solution of the problem (\ref{gr1}).

\section{Proof of the Theorem \protect\ref{1.1.1}}

To prove the existence of solution to (\ref{A1}) we consider the following
boundary value problem%
\begin{equation}
\begin{array}{ll}
-\Delta u+c(x)u^{-1}\left\vert \nabla u\right\vert ^{2}=a(x), & u>0\text{ in 
}B_{k},\text{ }u=0\text{ on }\partial B_{k},\text{ }%
\end{array}
\label{dd101}
\end{equation}%
where $B_{k}:=\{x\in \mathbb{R}^{N}\left\vert \left\vert x\right\vert
<k\right. \}$ is the ball of center $0$ and radius $k=1,2,..$. Put $\Omega
=B_{k}$ in \textit{Theorem \ref{idea1}. }Then the problem (\ref{dd101}) has
at least one solution $u_{k}\in C(\overline{B}_{k})\cap C^{2,\alpha }(B_{k})$%
, which satisfies%
\begin{equation}
\underline{u}_{2}\leq u_{k}\leq \overline{u}^{2}\text{ }in\text{ }B_{k},
\label{subsuper}
\end{equation}%
for $\underline{u}_{2}$ (resp. $\overline{u}^{2})$ the corresponding
functions from \textit{Theorem \ref{idea1} }when\textit{\ }$\Omega =B_{k}$%
\textit{.} In outside of $B_{k}$ we put $u_{k}=0$. The resulting function is
in $\mathbb{R}^{N}$. Now, we observe that 
\begin{equation}
w(r):=\int_{r}^{\infty }\xi ^{1-N}\int_{0}^{\xi }\sigma ^{N-1}\varphi
(\sigma )d\sigma d\xi ,\text{ }r:=\left\vert x\right\vert  \label{el2}
\end{equation}%
is the unique radial solution of the problem $-\Delta w=\varphi (\mid x\mid
) $ in $\mathbb{R}^{N}$, $w>0$ in $\mathbb{R}^{N}$, $w\overset{\left\vert
x\right\vert \rightarrow \infty }{\rightarrow }0$. \bigskip We prove that $w$
is bounded. Using integration by parts and L' Hôpital rule, we have%
\begin{eqnarray}
&&\int_{r}^{\infty }\xi ^{1-N}\int_{0}^{\xi }\sigma ^{N-1}\varphi (\sigma
)d\sigma d\xi =-\frac{1}{N-2}\int_{r}^{\infty }\frac{d}{d\xi }\left( \xi
^{2-N}\right) [\int_{0}^{\xi }\sigma ^{N-1}\varphi (\sigma )d\sigma ]d\xi 
\notag \\
&=&\frac{1}{N-2}\underset{R\rightarrow \infty }{\lim }\left\{
\int_{r}^{R}\xi \varphi (\xi )d\xi -R^{2-N}\int_{0}^{R}\sigma ^{N-1}\varphi
(\sigma )d\sigma +r^{2-N}\int_{0}^{r}\sigma ^{N-1}\varphi (\sigma )d\sigma
\right\}  \notag \\
&=&\frac{1}{N-2}\underset{R\rightarrow \infty }{\lim }\frac{%
R^{N-2}[\int_{r}^{R}\xi \varphi (\xi )d\xi +r^{2-N}\int_{0}^{r}\xi
^{N-1}\varphi (\xi )d\xi ]-\int_{0}^{R}\xi ^{N-1}\varphi (\xi )d\xi }{R^{N-2}%
}  \notag \\
&=&\frac{1}{N-2}\left[ \int_{r}^{\infty }\xi \varphi (\xi )d\xi
+r^{2-N}\int_{0}^{r}\xi ^{N-1}\varphi (\xi )d\xi \right] ,\text{ }R>r.
\label{c44}
\end{eqnarray}%
Now, by the second mean value theorem \ for integrals follows that there
exists $r_{1}\in (0,r)$ such that%
\begin{eqnarray}
\int_{0}^{r}\xi ^{N-1}\varphi (\xi )d\xi &=&\int_{0}^{r}\xi ^{N-2}\xi
\varphi (\xi )d\xi  \notag \\
&=&r^{N-2}\int_{r_{1}}^{r}\xi \varphi (\xi )d\xi \leq r^{N-2}\int_{0}^{r}\xi
\varphi (\xi )d\xi  \label{c55}
\end{eqnarray}%
for $N>2$. By\ (\ref{c44})-(\ref{c55}) we obtain $w(r)\leq K:=\frac{1}{N-2}%
\int_{0}^{\infty }\xi \varphi (\xi )d\xi $. We observe, in addition, that $w$
satisfies $-\Delta w(\left\vert x\right\vert )+c(x)w^{-1}(\left\vert
x\right\vert )\left\vert \nabla w(\left\vert x\right\vert )\right\vert
^{2}\geq a(x)$, $x\in \mathbb{R}^{N}$, $0<w\leq K$ and $w(r)\rightarrow 0$
as $r\rightarrow \infty $.

We prove that 
\begin{equation}
u_{k}\leq w(\left\vert x\right\vert )\text{, \ }x\in \mathbb{R}^{N}\text{, }%
k=1,2,3,...\text{ }  \label{r}
\end{equation}%
Since $w(\left\vert x\right\vert )>0$ in $\mathbb{R}^{N}$ and $u_{k}=0$ in $%
\mathbb{R}^{N}\backslash B_{k}$ it is enough to prove that $u_{k}\leq w$ in $%
B_{k},$ $k=1,2,3,...$ To prove this we observe that $w\in C^{2}\left( 
\overline{B}_{k}\right) $ and 
\begin{equation*}
\left\{ 
\begin{array}{llll}
-\Delta \lbrack w(x)-u_{k}(x)] & \geq & c(x)u_{k}^{-1}(x)\left\vert \nabla
u_{k}(x)\right\vert ^{2}-a(x)+a(x)\geq 0, & in\text{\ }B_{k}, \\ 
\text{ \ \ \ \ \ }w(x)-u_{k}(x) & > & 0, & on\text{ }\partial B_{k}.%
\end{array}%
\right.
\end{equation*}%
As a consequence of the maximum principle, \textit{Lemma \ref{maximum},} we
have that $u_{k}\leq w$ in $B_{k}$. So (\ref{r}) holds.

To finish the proof, use the standard convergence procedure (see \cite{CDC5}
or \cite{NS2}) and so $u_{k}$ has a subsequence, denoted again by $u_{k},$
such that $\displaystyle u_{k}\rightarrow u$ (pointwise) in $%
C_{loc}^{2,\alpha }({\mathbb{R}}^{N})$ and that $u$ is a solution for the
problem (\ref{gr1}) that vanishing at infinity.

In order to show (\ref{as22}), from the above arguments we have%
\begin{equation}
u\leq w\text{ in }\mathbb{R}^{N}\text{.}  \label{u22}
\end{equation}%
On the other hand, using (\ref{el2}) we have%
\begin{eqnarray*}
\lim_{\left\vert x\right\vert \rightarrow \infty }\frac{w(\left\vert
x\right\vert )}{\left\vert x\right\vert ^{2-\mu }} &=&\frac{1}{2-\mu }%
\lim_{\left\vert x\right\vert \rightarrow \infty }\frac{w^{\prime }(x)}{%
\left\vert x\right\vert ^{1-\mu }}=\frac{1}{\mu -2}\lim_{\left\vert
x\right\vert \rightarrow \infty }\left[ \int_{0}^{\left\vert x\right\vert
}\sigma ^{N-1}\varphi (\sigma )d\sigma /\left\vert x\right\vert ^{N-\mu }%
\right] \\
&=&\frac{1}{\mu -2}\lim_{\left\vert x\right\vert \rightarrow \infty
}\left\vert x\right\vert ^{\mu }\varphi (\left\vert x\right\vert )<\infty 
\text{.}
\end{eqnarray*}%
The above relation imply 
\begin{equation}
w(x)=O(\left\vert x\right\vert ^{2-\mu })\text{ }as\text{ }\left\vert
x\right\vert \rightarrow \infty .  \label{asimp}
\end{equation}%
Now, (\ref{as22}) follows from (\ref{asimp}) and (\ref{u22}). The proof of 
\textit{Theorem \ref{1.1.1}} is completed.

Affiliation 

Drago\c{s}-P\u{a}tru Covei

$^{1}$Constantin Brancusi University of Tg-Jiu and $^{2}$West University of
Timi\c{s}oara

$^{1}$Calea Eroilor, No 30, Tg-Jiu, Gorj, Romania and $^{2}$Bld. Pârvan, No.
4, Timi\c{s}oara, Timi\c{s}, Romania

e-mail: covdra@yahoo.com

\end{document}